\documentclass[12]{article}
\usepackage{amssymb,latexsym}
\usepackage{fullpage}
\usepackage{array}
\begin{document}
\title{K-theory of torus manifolds} 
\author{ V.Uma}
\date{} \maketitle 

\thispagestyle{empty} 


\def\theequation {\arabic{section}.\arabic{equation}}

\newcommand{\codim}{\mbox{{\rm codim}$\,$}}
\newcommand{\stab}{\mbox{{\rm stab}$\,$}}
\newcommand{\lr}{\mbox{$\longrightarrow$}}

\newcommand{\ch}{{\cal H}}
\newcommand{\cf}{{\cal F}}
\newcommand{\cd}{{\cal D}}

\newcommand{\ra}{\rightarrow}
\newcommand{\blr}{\Big \longrightarrow}
\newcommand{\da}{\Big \downarrow}
\newcommand{\ua}{\Big \uparrow}
\newcommand{\hra}{\mbox{\LARGE{$\hookrightarrow$}}}
\newcommand{\rt}{\mbox{\Large{$\rightarrowtail$}}}
\newcommand{\dua}{\begin{array}[t]{c}
\Big\uparrow \\ [-4mm]
\scriptscriptstyle \wedge \end{array}}

\newcommand{\be}{\begin{equation}}
\newcommand{\ee}{\end{equation}}

\newtheorem{guess}{Theorem}[section]
\newcommand{\bth}{\begin{guess}$\!\!\!${\bf .}~}
\newcommand{\eeth}{\end{guess}}
\renewcommand{\bar}{\overline}
\newtheorem{propo}[guess]{Proposition}
\newcommand{\bpropo}{\begin{propo}$\!\!\!${\bf .}~}
\newcommand{\epropo}{\end{propo}}

\newtheorem{lema}[guess]{Lemma}
\newcommand{\blem}{\begin{lema}$\!\!\!${\bf .}~}
\newcommand{\elem}{\end{lema}}

\newtheorem{defe}[guess]{Definition}
\newcommand{\bdefe}{\begin{defe}$\!\!\!${\bf .}~}
\newcommand{\edefe}{\end{defe}}

\newtheorem{coro}[guess]{Corollary}
\newcommand{\bcor}{\begin{coro}$\!\!\!${\bf .}~}
\newcommand{\ecor}{\end{coro}}

\newtheorem{rema}[guess]{Remark}
\newcommand{\brem}{\begin{rema}$\!\!\!${\bf .}~\rm}
\newcommand{\erem}{\end{rema}}

\newtheorem{exam}[guess]{Example}
\newcommand{\beg}{\begin{exam}$\!\!\!${\bf .}~\rm}
\newcommand{\eeg}{\end{exam}}

\newcommand{\ctext}[1]{\makebox(0,0){#1}}
\setlength{\unitlength}{0.1mm}
\newcommand{\cl}{{\cal L}}
\newcommand{\cp}{{\cal P}}
\newcommand{\ci}{{\cal I}}
\newcommand{\bz}{\mathbb{Z}}
\newcommand{\cs}{{\cal s}}
\newcommand{\cv}{{\cal V}}
\newcommand{\ce}{{\cal E}}
\newcommand{\ck}{{\cal K}}
\newcommand{\cR}{{\cal R}}
\newcommand{\cz}{{\cal Z}}
\newcommand{\cg}{{\cal G}}
\newcommand{\bq}{\mathbb{Q}}
\newcommand{\bt}{\mathbb{T}}
\newcommand{\bh}{\mathbb{H}}
\newcommand{\br}{\mathbb{R}}
\newcommand{\wt}{\widetilde}
\newcommand{\im}{{\rm Im}\,}
\newcommand{\bc}{\mathbb{C}}
\newcommand{\bp}{\mathbb{P}}
\newcommand{\spin}{{\rm Spin}\,}
\newcommand{\ds}{\displaystyle}
\newcommand{\tor}{{\rm Tor}\,}
\newcommand{\bff}{{\bf F}}
\newcommand{\bs}{\mathbb{S}}
\def\ns{\mathop{\lr}}
\def\nssup{\mathop{\lr\,sup}}
\def\nsinf{\mathop{\lr\,inf}}
\renewcommand{\phi}{\varphi}
\newcommand{\co}{{\cal O}}
\noindent
\begin{abstract} The {\it torus manifolds} have been defined and studied 
by Masuda and Panov (\cite {mp}) who in particular also describe its
cohomology ring structure.  In this note we shall describe the
topological $K$-ring of a class of torus manifolds (those for which
the orbit space under the action of the compact torus is a {\it
homology poytope} whose {\it nerve} is {\it shellable}) in terms of
generators and relations. Since these torus manifolds include the
class of quasi-toric manifolds this is a generalisation of our earlier
results (\cite {su}).
\end{abstract}

\section{Introduction}

We shall first briefly recall the notations and basic definitions from
\cite{mp}.

Let $M$ be a $2n$-dimensional closed connected orientable smooth
manifold with an effective smooth action of an $n$-dimensional torus
$T=(\bs^1)^n$ such that $M^{T}\neq \emptyset$. Since dim$(M)=2$ dim$(T)$
and $M$ is compact, the fixed point set $M^{T}$ is a finite set of
isolated points.

A closed connected codimension-two submanifold of $M$ is called {\it
characteristic} if it is the connected component of the set fixed
pointwise by a certain circle subgroup of $T$ and contains at least
one $T$-fixed point. Since $M$ is compact there are only finitely many
characteristic submanifolds. We denote them by $M_i
(i=1,\ldots,m)$. Note that each $M_i$ is orientable. We say that $M$
is {\it omnioriented} if an orientation is fixed for $M$ and for every
characteristic submanifold $M_i$. Further, $M$ is called a {\it torus
manifold} when it is omnioriented. 

Let $Q:=M/T$ denote the orbit space of $M$ and $\pi:M\ra Q$ the
quotient projection. We define the {\it facets} of $Q$ to be the orbit
spaces of the characteristic submanifolds:
$Q_i:=\pi(M_i),i=1,\ldots,m.$ Every facet is a closed connected subset
in $Q$ of codimension $1$.  We refer to a non-empty intersection of
$k$-facets as a codimension-$k$ {\it preface}, $k=1,\ldots, n$. Hence
a preface is the orbit space of some non-empty intersection
$M_{i_1}\cap\cdots\cap M_{i_k}$ of characteristic submanifolds. We
refer to the connected components of prefaces as {\it faces}. We also
regard $Q$ itself as a codimension-zero face; other faces are called
{\it proper faces}. A space $X$ is {\it acyclic} if $\wt{H}_i(X)=0$
for all $i$. We say that $Q$ is {\it face-acyclic} if all of its faces
(including $Q$ itself) are acyclic. We call $Q$ a {\it homology
polytope} if all its prefaces are acyclic (in particular,
connected). Note that $Q=M/T$ is a {\it homology polytope} if and only
if it is face-acyclic and all non-empty multiple intersections of
characteristic submanifolds $M_i$ are connected.

We say that a torus manifold $M$ is locally standard if every point in
$M$ has an invariant neighbourhood $U$ weakly equivariantly
diffeomorphic to an open subset $W\subset \bc^n$(invariant under the
standard $T$-action on $\bc^n$). The latter means that there is an
automorphism $\psi:T\ra T$ and a diffeomorphism $f:U\ra W$ such that
$f(ty)=\psi(t)f(y)$ for all $t\in T,y\in U$.

Any point in the orbit space $Q$ of a locally standard torus manifold
$M$ has a neighbourhood diffeomorphic to an open subset in the
positive cone $$\br^n_{+}=\{(y_1,\ldots,y_n)\in \br:y_i\geq 0,i=1\ldots,n\}.$$

Moreover, this local diffeomorphism preserves the face structures in
$Q$ and $\br^n_{+}$(that is, a point from a codimension-k face of $Q$
is mapped to a point with at least $k$ zero coordinates). By the
definition, this identifies $Q$ as a {\it manifold with corners}. In
particular $Q$ is a manifold with boundary $\partial Q=\cup_i Q_i$. Let
$K$ denote the $nerve$ of the covering of $\partial Q$ by the
facets. Thus $K$ is an $(n-1)$-dimensional simplicial complex on
$m$-vertices. The $(k-1)$-dimensional simplices of $K$ are in one-one
correspondence with the codimension-$k$ prefaces of $Q$.

We assume that a torus manifold $M$ is locally standard. Then the
orbit space $Q$ is a manifold with corners. The facets of $Q$ are the
quotient images $Q_i$ of characteristic submanifolds $M_i (i=1,\ldots
m)$. Let $\Lambda:\{1,\ldots,m\}\ra H_2(BT)=Hom(\bs^1,T)\simeq \bz^n$
be a map sending $i$ to $a_i$ where the circle subgroup determined by
$a_i$, that is $a_i(\bs^1)$, is the one which fixes $M_i$ (see
Prop. 2.5 and \S3.2 of \cite{mp}). Further, the {\it characteristic}
map $\Lambda$ satisfies the following {\it non-singular condition}:

If $Q_{i_1}\cap\cdots \cap Q_{i_k}$ is non-empty, then
$\Lambda(i_1),\ldots,\Lambda(i_k)$ span a $k$-dimensional unimodular
subspace (i.e extend to a $\bz$-basis) of $Hom(\bs^1,T)\simeq \bz^n$.

The data $(Q,\Lambda)$ determines the torus manifold $M$ if the orbit
quotient $Q$ of $M$ satisfies $H^2(Q)=0$ (see Lemma 3.6 of \cite{mp}
and Prop. 1.8 of \cite{dj}).

Let $Q$ be a homology polytope (or even a simple convex polytope) with
$m$ facets $Q_1,\ldots, Q_m$. Let ${\bf k}$ be a commutative ground
ring with unit. Then the Stanley-Reisner face ring of its nerve $K$
can be identified with the ring
$${\bf k}[Q]={\bf k}[v_{Q_1},\ldots,v_{Q_m}]/(v_{Q_{i_1}}\cdots
v_{Q_{i_k}}~ if~ Q_{i_1}\cap\cdots Q_{i_k}=\emptyset)$$ called the
{\it face ring} of $Q$ (see \S4 of \cite{mp}).

\section{Main Theorem}

This section is devoted to proving our main result Theorem \ref{main}.

\bpropo\label{lbs}There exists a complex line bundle $\cl_j$ on $M$
admitting a section $s_j$ whose zero locus is the characteristic
submanifold $M_j$ for $1\leq j\leq m$.  \epropo {\bf Proof:} Let
${\nu}_{j}$ denote the normal bundle of $M_j$ in $M$ and let ${\sf p}:
{\nu_j}: \rightarrow M_j$ be the canonical projection and let $E(\nu_j)$
denote the total space of $\nu_j$. The rank $2$ real vector bundle
$\nu_j$ on $M_j$ admits a Riemannian metric (since $M_j$ is compact)
and in fact its structure group can be reduced to $O(2)$. Fixing
orientations on $M$ and $M_j$ determines a canonical orientation for
every $\nu_j$ for $1\leq j\leq m$. Therefore the normal bundle $\nu_j$
admits reduction of structure group to $SO(2)$. We can identify
$SO(2)$ with $S^1$ so that the principal $SO(2)$ bundle associated to
$\nu_j$ is in fact an $S^1$ bundle and the complex line bundle
associated to it by the standard action of $S^1$ on $\bc$ has $\nu_j$
as its underlying real vector bundle.

By the tubular neighbourhood theorem, $E(\nu_j)$ is diffeomorphic to a
tubular neighbourhood $D_j$ of $M_j$ and the diffeomophism maps the
image of the zero section of $\nu_j$ onto $M_j$. Further, {\it the
total space} of the principal $S^1$ bundle can be identified with
$\partial(D_j)$, the boundary of the tubular neighbourhood. Let
${\sf p}^*(\nu_j)$ be the pull back of $\nu_j$ to $D_j$. Since $\nu_j$ is
associated to the principal $S^1$ bundle, its pull back to
${\partial(D_j)}$ (its total space) is trivial and since
$\partial(D_j)$ is a deformation retract of $D_j- M_j$, by the
homotopy property of vector bundles it follows that the restriction of
${\sf p}^*(\nu_j)$ to $D_j- M_j$ is {\em trivial}.

The vector bundle ${\sf p}^*(\nu_j)$ is endowed with a section $\sigma_j$
namely the diagonal section whose zero locus $Z(\sigma_j)= M_j$.

Let $\epsilon$ be the trivial complex line bundle on
$M-int(D_j)$. Note that ${\sf p}^*(\nu_j)$ and $\epsilon$ agree on a
neighbourhood of $ D_j\cap (M-int(D_j))=\partial (D_j)$. Thus we can
construct a complex line bundle ${\cal L}_j $ on the whole of $M$
which agrees with ${\sf p}^*(\nu_j)$ on $D_j$ and with $\epsilon$ on
$M-int(D_j)$ (see \cite{k}). Further, the section $\sigma_j$ of
$\nu_j$ extends to give a section $s_j$ for $\cl_j$ whose zero locus $
Z(s_j)= M_j$.\hfill $\Box$

We now recall the notion of a {\it shellable simplicial complex} (see
Def 2.1, page 79 of \cite{st}).

A simplicial complex $\Delta$ is said to be {\it pure} if each of its
facets (or maximal face) has the same dimension. We say that a pure
simplicial complex $\Delta$ is {\it shellable} if its facets can be
ordered $F_1,\ldots,F_s$ such that the following condition holds: Let
$\Delta_j$ be the subcomplex generated by $F_1,\ldots,F_j$, i.e:
$$\Delta_j=2^{F_1}\cup\cdots\cup 2^{F_j}$$ where $2^{F}=\{G:G\subseteq
F\}$. Then we require that for all $1\leq i\leq s$ the set of faces of
$\Delta_i$ which do not belong to $\Delta_{i-1}$ has a unique minimal
element (with respect to inclusion). (When $i=1$, we have
$\Delta_0=\emptyset$ and $\Delta_1=2^{F_1}$, so $\Delta_1-\Delta_0$
has the unique minimal element $\emptyset$.) The linear order
$F_1,\ldots,F_s$ is called a {\it shelling order} or a {\it shelling}
of $\Delta$. Given a shelling $F_1,\ldots,F_s$ of $\Delta$, we define
the {\it restriction} $r(F_i)$ of $F_i$ to be the unique minimal
element of $\Delta_i- \Delta_{i-1}$.

{\it Henceforth we assume that $M$ is a torus manifold with orbit
space a homology polytope $Q$ whose nerve $K$ is a shellable
simplicial complex.}

Let $d$ be the number of vertices of $Q$ so that the simplicial
complex $K$ has $d$ {\it facets} (or maximal dimensional faces).  Let
$F_1,\ldots, F_d$ be a shelling of $K$ (see Def 2.1, page 79 of
\cite{st}) and let $r(F_i)$ denote the {\it restriction} of
$F_i$. Thus (by immediate consequence of the definition of a shelling)
we have a disjoint union:
$$K=[r(F_1),F_1]\sqcup\cdots\sqcup[r(F_d),F_d].$$ 
Let $S_1,\ldots
S_d$ denote the vertices of the polytope which correspond respectively
to $F_1,\ldots,F_d$.  Further, let $T_i$ be the face of $Q$
corresponding to the face $r(F_i)$ of $K$ for $1\leq i\leq d$.  Thus
it follows that every face of $Q$ belongs to $[S_i,T_i]$ for a unique
$1\leq i\leq d$ (where $[S_i,T_i]$ stands for the faces of $Q$ which
contain the vertex $S_i$ and lie on the face $T_i$). We isolate this
property as follows: For every face $Q_I=Q_{j_1}\cap\cdots\cap
Q_{j_k}$ of $Q$ {where $I=\{j_1,\ldots,j_k\}$ there is a unique $1\leq
i\leq d$ such that:
$$Q_I\in[S_i,T_i]\eqno{(*)}$$

Let $\tau_i:=dim(T_i)$.  Further, if $\widehat{T_i}$ denotes the
subset of $T_i$ obtained by deleting all faces of $T_i$ not containing
$S_i$. Since $Q$ is a manifold with corners, $\widehat{T_i}$ is
identified with $\br_{+}^{\tau_i}$.  Then we note that
$\pi^{-1}(\widehat{T_i})$ identified with $\bc^{\tau_i}$ for $1\leq
i\leq d$ give a cellular decomposition of $M$ where $\pi:M\ra Q$ is
the quotient projection (see Construction 5.15 on page 66 of \cite{bp}
and \S3 of \cite{dj}).  Hence we can summarise as follows:
\blem{\label{shell}} Let $M$ be a torus manifold with orbit space $Q$
a homology polytope. Let $K$ be the nerve of $Q$. If $K$ is shellable
then the shelling gives a perfect cellular decomposition of $M$ with
cells only in even dimensions.  \elem

\bth \label{main} Let $Q_1,\ldots, Q_m$ denote the facets of
$Q$. Let $a_i$ denote the element $\Lambda(i)$ in $H_{2}(BT)$, where
$\Lambda$ is the characteristic homomorphism.  Consider the polynomial
algebra $\bz[v_{Q_1},\cdots, v_{Q_m}]$.  We denote by $I$ the ideal
generated by the following two types of elements
 $$v_{Q_{j_1}}\cdots v_{Q_{j_k}}, ~1\leq j_p\leq m, \eqno{(i)} $$
 where
 $\cap_{i=1}^k Q_{j_i}= \emptyset$ in $Q$, and the elements
\noindent
$$\prod_{j, \langle t,a_j\rangle>0} (1-v_{Q_j})^{\langle
t,a_{j}\rangle}- \prod_{j, \langle
t,a_{j}\rangle<0}(1-v_{Q_j})^{-\langle t,a_j\rangle}\eqno {(ii)}$$ for
$t\in H^2(BT)$.  Let $\cR=\bz[v_{Q_1},\cdots, v_{Q_m}]/I$ and let
$K^{*}(M)$ denote the topological K-ring of $X$. Then the map
$\psi:\cR\rightarrow K^{*}(M)$ sending $v_{Q_j}$ to $[\cl_{j}]-1$ is a
ring isomorphism.  \eeth

We now state the following lemma which is used for proving Theorem
\ref{main}.  \blem\label{kring} The monomials $v_{T_i}, 1\leq i\leq d$
span $\cR$ as a $\bz$-module.  \elem {\bf Proof:} The proof of this
lemma is exactly as of Prop. 2.1 of \cite{su} (also see Lemma 2.2 of
\cite{su1}). Thanks to the key observation ${(*)}$, the arguments work
analogously in this setting too (the setting of a {\it torus manifold}
with quotient a {\it homology polytope} whose {\it nerve} is {\it
shellable}). We omit the details.\hfill $\Box$

{\bf Proof of Theorem \ref{main}:} By Lemma \ref{lbs}, there exists a
complex line bundle $\cl_j$ on $M$ with section $s_j$ whose zero locus
$Z(s_j)=M_j$. Thus its first chern class, $c_1(\cl_j)=[M_j]$ for
$1\leq j\leq m$, where $[M_j]$ denotes the fundamental class of $M_j$
in $H^2(M;\bz)$. Further, $M$ has a cellular decomposition with cells
only in even dimensions (see Lemma \ref{shell}). Now, by Corollary 6.8
of \cite{mp}, $H^*(M;\bz)$ is generated by $c_1(\cl_1),\ldots
c_1(\cl_m)\in H^2(X;\bz)$. Hence by Theorem 4.1 of \cite{su1}, it
follows that $K^*(M)=K^0(M)$ is generated by $[\cl_1],\ldots,
[\cl_m]\in K^0(M)$.

Let $\cl_t:=\prod_{j=1}^m \cl_j^{\langle t,a_j\rangle}$. Since
$c_1(\cl_t)=\sum _{j=1}^m \langle t,a_j\rangle [M_j]=0$ it follows
that $\cl_t$ is a trivial line bundle for every $t\in H^2(BT)$.

Let $Q_{j_1}\cap\cdots\cap Q_{j_k}=\emptyset$ in $Q$. Then we have $
M_{j_1}\cap\cdots\cap M_{j_k}=\emptyset$ in $M$.  Therefore the vector
bundle $\cv=\cl_{j_1}\oplus\cdots\oplus\cl_{j_k}$ admits a section
$s=(s_{j_1},\ldots, s_{j_k})$ which is nowhere vanishing. Hence
applying the $\gamma^k$ operation in $K(M)$ we obtain
$\gamma^k(\oplus_{p=1}^k \cl_{j_p}-k)=\gamma^k(\oplus_{p=1}^k
(\cl_{j_p}-1))= \prod_{p=1}^k ([\cl_{j_p}]-1)$. Since $\cv$ has
geometric dimension at most $k-1$ we have: $\prod_{p=1}^k
([\cl_{j_p}]-1)=0$.

By the above arguments it follows that the map $\psi:\cR\rightarrow
K^*(X)$ which sends $v_{Q_j}$ to $[\cl_j]-1$ is well defined and
surjective.

Since $M$ has a cell decomposition with cells only in even dimensions
 it follows that $K^*(M)=K^0(M)$ is free abelian of rank $d$ which is
 the number of even dimensional cells (see \cite{ah}). Further, by
 Lemma \ref{kring} we know that $\cR$ is generated by $d$ elements
 $v_{T_1},\ldots, v_{T_d}$. Hence it follows that the map $\psi$ is a
 ring isomorphism. \hfill $\Box$

\noindent
{\it Acknowledgement:} I thank P. Sankaran for helpful discussions.


\noindent
Department of Mathematics\\
I.I.T Madras,\\
Chennai-600 036\\
INDIA.\\
E-mail:{\tt vuma@iitm.ac.in}\\
{\tt uma@cmi.ac.in}}

\end{document}